\documentclass{tran-l}
 \theoremstyle{definition}
 
 \theoremstyle{remark}
 
 \numberwithin{equation}{section}

\begin{document}

\title[Euler Type Generalization of Wilson's Theorem]
{Euler Type Generalization of Wilson's Theorem}

\author{Mehdi Hassani\\ Mahmoud Momeni-Pour}

\address{Institute for Advanced Studies in Basic Sciences\\
P.O. Box 45195-1159\\
Zanjan, Iran.}

\email{mmhassany@srttu.edu, m$_{-}$momeni@iasbs.ac.ir}

\thanks{}



\date{}

\dedicatory{}

\commby{}
%%% ----------------------------------------------------------------------
\begin{abstract}
In this short note, we introduce an Euler analogue of Wilson's
theorem; $a_1a_2\cdots a_{\phi(n)}\equiv (-1)^{\phi(n)+1}~({\rm
mod}~n)$ say, where ${\rm gcd}(a_i,n)=1$.

\end{abstract}

\maketitle

%\tableofcontents

%\section{Introduction and Results}
Recently, some generalizations of Wilson's theorem \cite{apos};
$(p-1)!\equiv-1~~({\rm mod}~p)$, which $p$ is a prime number, has
been taken for the nonzero elements of a finite field \cite{has}.
In this short note, we get it for the elements of a finite
multiplicative subgroup of a commutative ring with identity. We
start with $\mathbb{Z}_n=\{0,1,2,\cdots,n-1\}$, and
$$
U(\mathbb{Z}_n)=\{a\in\mathbb{Z}_n: {\rm
gcd}(a,n)=1\}=\{a_1,a_2,\cdots,a_{\phi(n)}\}.
$$
Since $U(\mathbb{Z}_n)$ is a multiplicative group, by Lagrange's
theorem, if $a\in U(\mathbb{Z}_n)$ then $o(a)|\phi(n)$ and so
$a^{\phi(n)}=1$. In the other hand,
$x^{\phi(n)}-1\in\mathbb{Z}[x]$ and $\mathbb{Z}$ is an integral
domain, so the elements of $U(\mathbb{Z}_n)$ are actually the
roots of $x^{\phi(n)}-1$, that is
$\mathcal{Z}(x^{\phi(n)}-1)=U(\mathbb{Z}_n)$. Thus, we obtain
$$
x^{\phi(n)}-1=\prod_{i=1}^{\phi(n)}(x-a_i).
$$
Considering elementary symmetric functions \cite{pat};
\begin{eqnarray*}
s_k=s_k(a_1,a_2,\cdots,a_{\phi(n)})=\sum_{1\leq i_1< i_2<\cdots<
i_k\leq \phi(n)}a_{i_1}a_{i_2}\cdots a_{i_k},
\end{eqnarray*}
we have
\begin{eqnarray*}
x^{\phi(n)}-1&=&(x-a_1)(x-a_2)\cdots(x-a_{\phi(n)})\\&=&x^{\phi(n)}-s_1
x^{\phi(n)-1}+s_2 x^{\phi(n)-2}+\cdots+(-1)^{\phi(n)}s_{\phi(n)},
\end{eqnarray*}
or the following identity
\begin{equation*}
\sum_{k=1}^{\phi(n)}(-1)^ks_kx^{\phi(n)-k}+1=0,
\end{equation*}
and comparing coefficients, we obtain
$$
s_1=0, s_2=0,\cdots, s_{\phi(n)-1}=0~{\rm
and}~s_{\phi(n)}=(-1)^{\phi(n)+1},
$$
which we can state all of them together as follows
\begin{equation*}
s_k=\left\lfloor\frac{k}{\phi(n)}\right\rfloor(-1)^{\phi(n)+1}\hspace{10mm}(k=1,2,\cdots,\phi(n)).
\end{equation*}
This relation is an Euler type generalization of the Wilson's
theorem. Specially, considering it for $k=\phi(n)$, we have
$$
\prod_{a\in U(\mathbb{Z}_n)}a\equiv
(-1)^{\phi(n)+1}\hspace{10mm}({\rm mod}~n).
$$
Letting $n$ an odd prime in above, reproves Wilson's theorem.
Finally, we mention that, more generally if $R$ is a commutative
ring with identity and $H$ is a finite multiplicative subgroup of
it, then using division algorithm \cite{hun}, similarly we obtain
$$
s_k(H)=\left\lfloor\frac{k}{|H|}\right\rfloor(-1)^{|H|+1}\hspace{10mm}(k=1,2,\cdots,|H|),
$$
where $s_k(H)$ is $k$-th elementary symmetric function of the
elements of $H$. Specially, putting $k=|H|$, we obtain the
following generalization of the Wilson's theorem
$$
\prod_{a\in H}a=(-1)^{|H|+1}.
$$

% ------------------------------------------------------------------------

\end{document}